\documentclass{article}

\usepackage{amsmath}\usepackage{amsfonts}\usepackage{amssymb}

\newtheorem{theorem}{Theorem}[section]
\newtheorem{prop}[theorem]{Proposition}
\newcommand{\qedwhite}{\hfill \ensuremath{\Box}}
\date{}
\title{On the Jacobi formula for Bivariate\\ Pad\'e Approximants of Rectangular Type}
\author{Gareth Hegarty}
\begin{document}

\maketitle
\begin{abstract}
In this paper a recursive algorithm is presented for evaluating multivariate
Pad\'e approximants (of the rectangular type described in the work of Lutterodt) 
which is analogous to the Jacobi formula for univariate Pad\'e approximants. This 
algorithm is then applied to a (singular) {\em Riccati differential equation} to generate fast and accurate
approximate solutions.
\end{abstract}


\section{Introduction}\label{s1}
Consider a bivariate power series for a function which vanishes at the origin:
\begin{equation}
f(x,y)=\sum^\infty_{n,m=0}c_{n,m}x^ny^m;\ c_{0,0}=0
\label{biseries}
\end{equation}
The aim of this paper is to develop recursive algorithms for evaluating diagonal Pad\'e approximants
based on this series. To set the stage for this, we first consider the special cases of the {\em univariate}
diagonal Pad\'e approximants in $x$ and $y$.
Thus, in $x$ we have the (diagonal) $[n/n]$-Pad\'e approximant for $f(x,0)$ as a ratio of polynomials of degree $n$:
$$f_{n,0}(x)=\frac{A_{n,0}(x)}{B_{n,0}(x)}=\frac{a^{(n)}_{1,0}x+...+a^{(n)}_{n,0}x^n}
{1+b^{(n)}_{1,0}x+...+b^{(n)}_{n,0}x^n}$$
which agrees with
the series for $f(x,0)$ up to order $x^{2n}$, i.e.
\begin{equation}
E_{n,0}=A_{n,0}-B_{n,0}C_{2n,0}=\sum^n_{k=1}e^{(n)}_{2n+k}x^{2n+k};\ C_{2n,0}=\sum^{2n}_{k=0}c_{k,0}x^k
\label{v1}
\end{equation}
One of the most interesting aspects of the diagonal Pad\'e approximants is the recursion
relation between successive numerators and denominators known as the {\em Jacobi formula}
(described in Fair, \cite{fair}), i.e.
starting with $A_{-1,0}=E_{-1,0}=-x^{-1}$, $B_{-1,0}=0$ and $A_{0,0}=E_{0,0}=0$, $B_{0,0}=1$, the
subsequent iterates satisfy:
\begin{eqnarray}
A_{n,0}=(1+\beta_n x)A_{n-1,0}+\alpha_n x^2 A_{n-2,0},\label{jaca}\\
B_{n,0}=(1+\beta_n x)B_{n-1,0}+\alpha_n x^2 B_{n-2,0},\label{jacb}
\end{eqnarray}
for $n>0$. The case $n=1$ is obvious, i.e. $\alpha_1=-c_{1,0}$, $\beta_1=c_{1,0}$. For $n>1$ we subtitute these 
equations back into (\ref{v1}) to give:
\begin{eqnarray}
E_{n,0}&=&(1+\beta_n x)\{E_{n-1,0}-B_{n-1,0}(c_{2n-1,0}x^{2n-1}+c_{2n,0}x^{2n})\}\label{En0}\\
&&+\alpha_n x^2 \{E_{n-2,0}
-B_{n-2,0}(c_{2n-3,0}x^{2n-3}+c_{2n-2,0}x^{2n-2}\nonumber\\
&&+c_{2n-1,0}x^{2n-1}+c_{2n,0}x^{2n})\}\nonumber
\end{eqnarray}
The coefficients $\alpha_n,\beta_n$ can be found from this,
i.e. $\alpha_n$ comes from the coefficient of $x^{2n-1}$ and $\beta_n$ from the coefficient of $x^{2n}$:
\begin{equation}
\alpha_n=-\big\{e^{(n-1)}_{2(n-1)+1}-b^{(n-1)}_{0,0}c_{2n-1,0}\big\}\big/
\big\{e^{(n-2)}_{2(n-2)+1}-b^{(n-2)}_{0,0}c_{2n-3,0}\big\}
\label{alpn}
\end{equation}
\begin{eqnarray}
\beta_n&=&\big\{b^{(n-1)}_{1,0}c_{2n-1,0}+b^{(n-1)}_{0,0}c_{2n,0}-e^{(n-1)}_{2(n-1)+2}\label{betn}\\
&&+\alpha_n\big[
b^{(n-2)}_{0,0}c_{2n-2,0}+b^{(n-2)}_{1,0}c_{2n-3,0}\nonumber\\
&&-e^{(n-2)}_{2(n-2)+2}\big]\big\}\big/
\big\{e^{(n-1)}_{2(n-1)+1}-b^{(n-1)}_{0,0}c_{2n-1,0}\big\}\nonumber
\end{eqnarray}
Similarly, in $y$ we have
the $[m/m]$-Pad\'e approximant for $f(0,y)$, i.e. $f_{0,m}(y)=A^{(m)}_{0,m}(y)/B^{(m)}_{0,m}(y)$,
following similar recursion relations to (\ref{jaca}), (\ref{jacb}) with corresponding parameters denoted
$\beta^*_m$, $\alpha^*_m$.

Acording to Cuyt, \cite{cuyt3}, there are several different approaches to defining (bivariate)
Pad\'e approximants for a bivariate power series such as (\ref{biseries}).
Here we follow the {\em equation lattice} approach as it is simple to relate to the univariate case
and is appropriate for the applications we consider. This involves choosing index sets in ${\mathbb N}^2_0$
(where ${\mathbb N}_0=\{0\}\cup {\mathbb N}$)
for the numerator, denominator, and equations such that the coefficients are well-defined
(see Levin, \cite{levin} for an outline of the
general framework).
Even within this approach, there are many ways of choosing such index sets (including the
diagonal approximants of Chisolm, \cite{chisolm}) but here we restrict our attention to a couple of definitions
of the {\em rectangular} type found in Lutterodt, \cite{lutt}.
Thus, for non-negative integers $n,m$, the {\em left}-$(n,m)$ Pad\'e approximant is formally defined as 
$f^L_{n,m}(x,y)=[N_{n,m}/D_{n,m}]_{Q_{n,m}}$ where
$$N_{n,m}=D_{n,m}=\{(i,j):i\leq n, j\leq m\}$$
$$Q_{n,m}= \{(i,j):\leq 0\leq i\leq 2n, j\leq m\}\cup \{(0,j): m+1\leq  j\leq 2m\}$$
In practical terms, this means that it may be written as a ratio of the form:
\begin{equation}
f^L_{n,m}(x,y)=\frac{A^{(n,m)}_{n,0}+ A^{(n,m)}_{n,1}y+...+A^{(n,m)}_{n,m}y^m}{
B^{(n,m)}_{n,0}+B^{(n,m)}_{n,1}y+...+B^{(n,m)}_{n,m}y^m}
\label{lpade}
\end{equation}
where $A^{(n,m)}_{n,p}$ and $B^{(n,m)}_{n,p}$ are polynomials in $x$ of degree $n$ for $p=0,...,m$, i.e. 
$A^{(n,m)}_{n,p}=a^{(n,m)}_{0,p}+...+a^{(n,m)}_{n,p}x^n$, 
$B^{(n,m)}_{n,p}=b^{(n,m)}_{0,p}+...+b^{(n,m)}_{n,p}x^n$ satisfying the equations
\begin{equation}
A^{(n,m)}_{n,p}-B^{(n,m)}_{n,0}C_{2n,p}-...-B^{(n,m)}_{n,p}C_{2n,0}= O(x^{2n+1});\ C_{2n,k}=\sum^{2n}_{i=0}c_{i,k}x^i
\label{v2}
\end{equation}
These equations are coupled with the conditions:
\begin{equation}
a^{(n,m)}_{0,p}=a^{(m)}_{0,p},\ b^{(n,m)}_{0,p}=b^{(m)}_{0,p},\, p=0,...,m
\label{conditions}
\end{equation}
where
$a^{(m)}_{0,p}$, $b^{(m)}_{0,p}$ are the $y^p$ coefficients in $A^{(m)}_{0,m}$, $B^{(m)}_{0,m}$
respectively.
Note that we can also define a corresponding {\em right}-$(n,m)$ Pad\'e approximant using symmetry, i.e.
\begin{equation}
f^R_{n,m}(x,y)=(f^T)^L_{m,n}(y,x);\ \ f^T(x,y)=f(y,x)
\label{rpade}
\end{equation}
Both $f^L_{n,m}$ and $f^R_{n,m}$ are special cases of the rectangular scheme found in Lutterodt, \cite{lutt},
and are diagonal in the sense that the order of dependence
on $x$ and $y$ in numerator and denominator are identical (but different to the diagonal Chisolm
approximants, \cite{chisolm}).
This is a flexible scheme because not only does it give two approximants to consider for each $n$, $m$
but also allows one to vary the order of dependence on $x$ and $y$ independently.

\section{Recursive Algorithm}
\label{s2}
Here is presented a recursive algorithm for generating the {\em left}-$(n,m)$ Pad\'e approximants
$f^{L}_{n,m}(x,y)$ over both $n$ and $m$. 
This will also generate the {\em right}-$(n,m)$ Pad\'e approximants
$f^{R}_{n,m}(x,y)$ due to the symmetry of the definition in (\ref{rpade}). 
In the case $m=0$, $f^{L}_{n,0}(x,y)=f_{n,0}(x)$ is just the $[n/n]$-Pad\'e approximant
of $f(x,0)$, to which the univariate Jacobi formulas apply. 
To facilitate the induction over $m$, we first generalise the definition above to allow the {\em seeds} of the
Pad\'e approximant $b_{0,1},...,b_{0,m}$ (for $m>0$) to be scalar unknowns (with $b_{0,0}=1$).
Inductively assuming that the previous polynomials 
$B_{n,p}=B_{n,p}[b_{0,1},...,b_{0,p}]$ 
(of degree $n$) 
for $p=0,...,m-1$ have been evaluated,
we define polynomials
$$A_{n,m}=A_{n,m}[b_{0,1},...,b_{0,m}]=a^{(n)}_{0,m}+a^{(n)}_{1,m}x+...+a^{(n)}_{n,m}x^n,$$
$$B_{n,m}=B_{n,m}[b_{0,1},...,b_{0,m}]=b^{(n)}_{0,m}+b^{(n)}_{1,m}x+...+b^{(n)}_{n,m}x^n,$$ 
both of degree $n$
(the latter with constant component $b^{(n)}_{0,m}=b_{0,m}$) 
satisfying the equation:
\begin{equation}
A_{n,m}-B_{n,0}C_{2n,m}-...-B_{n,m-1}C_{2n,1}-B_{n,m}C_{2n,0}= O(x^{2n+1}),
\label{v8}
\end{equation}
i.e. using the $x^{n+1},...,
x^{2n}$ equations to solve for the coefficients in $B_{n,m}$ in terms of the unknown parameter
$b_{0,m}$, which can then be fed into the $x^0,...,x^n$ equations for the coefficients in $A_{n,m}$.
It is obvious by construction that using the values $b_{0,p}=b^{(m)}_{0,p}$ for $p=1,...,m$ (the $y^p$
coefficients in $B^{(m)}_{0,m}$)
one recovers $f^L_{n,m}(x,y)$, i.e.
\begin{equation}
A_{n,p}[b^{(m)}_{0,1},...,b^{(m)}_{0,p}]=A^{(n,m)}_{n,p},\ 
B_{n,p}[b^{(m)}_{0,1},...,b^{(m)}_{0,p}]=B^{(n,m)}_{n,p},
\label{evalAB}
\end{equation}
for $p=1,...,m$.

The induction over $n$ requires a different approach that is consistent with the univariate Jacobi formulas in 
(\ref{jaca}), (\ref{jacb}). 
Starting with the initial values $A_{-1,m}=B_{-1,m}=0$ and
$A_{0,m}=a_{0,m}=c_{0,m}+b_{0,1}c_{0,m-1}+...+b_{0,m-1}c_{0,1}$,
$B_{0,m}=b_{0,m}$ we generate subsequent interates according to the following proposition: 
\begin{prop}
The polynomials $A_{n,m}$, $B_{n,m}$ for $n>0$ satisfy recursive relations of the form:
\begin{eqnarray}
A_{n,m} &=& (1+\beta_n x)A_{n-1,m}+\alpha_nx^2A_{n-2,m}+\check A_{n,m}
\label{Anm}\\
B_{n,m} &=& (1+\beta_n x)B_{n-1,m}+\alpha_nx^2B_{n-2,m}+\check B_{n,m}
\label{Bnm}
\end{eqnarray}
where
\begin{equation}
\check A_{n,m}
=\sum^{m-1}_{p=0}x\check\beta^{n,m}_pA_{n-1,p}
+x^2\check\alpha^{n,m}_pA_{n-2,p},
\label{checkAnm}
\end{equation}
\begin{equation}
\check B_{n,m}=\sum^{m-1}_{p=0}x\check\beta^{n,m}_pB_{n-1,p}
+x^2\check\alpha^{n,m}_pB_{n-2,p}
\label{checkBnm}
\end{equation}
The coefficients $\check\beta^{n,m}_0$, 
$\check\alpha^{n,m}_0$ are given in (\ref{checkbetnm}), (\ref{checkalphanm}) while $\check\beta^{n,m}_p$, 
$\check\alpha^{n,m}_p$ for $p=1,...,m-1$ are just copies of previous coefficients
in the sense that they satisfy the identities:
\begin{equation}
\check\beta^{n,m}_p=\check\beta^{n,m-s}_{p-s},\ 
\check\alpha^{n,m}_p=\check\alpha^{n,m-s}_{p-s};\ s=1,...,p,\ p=1,...,m-1,
\label{recidab}
\end{equation}
\end{prop} 
{\noindent\bf Proof.}
To justify these rules we consider the error polynomial:
$$
E_{n,m}=A_{n,m}-B_{n,0}C_{2n,m}-...-B_{n,m}C_{2n,0},
$$
and aim to show that if $A_{n,m}$, $B_{n,m}$ are evaluated according to the
proposition then $E_{n,m}=O(x^{2n+1})$, or more specifically
\begin{equation}
E_{n,m}=e^{n,m}_{2n+1}x^{2n+1}+...+e^{n,m}_{3n}x^{3n}
\label{Enm}
\end{equation}
Expanding $E_{n,m}$ using the rules in (\ref{Anm}), (\ref{Bnm}) and using the starting values of $E_{-1,m}=E_{0,m}=0$ we get:
$$
E_{n,m}=
(1+\beta_nx)\{E_{n-1,m}-
\sum^m_{p=0}B_{n-1,p}(c_{2n-1,m-p}x^{2n-1}
+c_{2n,m-p}x^{2n})\}
$$
$$
+\alpha_nx^2\{E_{n-2,m}-\sum^m_{p=0}B_{n-2,p}(c_{2n-3,m-p}x^{2n-3}
+c_{2n-2,m-p}x^{2n-2}
$$
$$+c_{2n-1,m-p}x^{2n-1}
+c_{2n,m-p}x^{2n})\}
+\check A_{n,m}-\check B_{n,1}C_{2n,m-1}-...--\check B_{n,m}C_{2n,0}
$$
This can be rewritten in the form:
$$
E_{n,m}=
-(\hat \sigma^{(m)}_{2n-1}x^{2n-1}+\hat \sigma^{(m)}_{2n}x^{2n})+
\hat\tau^{n,m}_{2n+1}x^{2n+1}+...+\hat\tau^{n,m}_{3n}x^{3n}
$$
$$
+\check A_{n,m}-\check B_{n,1}C_{2n,m-1}-...--\check B_{n,m}C_{2n,0}
$$
where
$$
\hat \sigma^{(m)}_{2n-1}=-e^{n-1,m}_{2n-1}+\sum^{m}_{p=0}(b^{(n-1)}_{0,p}c_{2n-1,m-p}+\alpha_nb^{(n-2)}_{0,p}c_{2n-3,m-p})
$$
$$
\hat \sigma^{(m)}_{2n}=-(e^{n-1,m}_{2n}+\beta_n e^{n-1,m}_{2n-1}+\alpha_n e^{n-2,m}_{2n-2})
+\sum^{m}_{p=0}\{b^{(n-1)}_{0,p}c_{2n,m-p}
$$
$$
+(b^{(n-1)}_{1,p}
+\beta_nb^{(n-1)}_{0,p})c_{2n-1,m-p}
+\alpha_n(b^{(n-2)}_{0,p}c_{2n-2,m-p}+b^{(n-2)}_{1,p}c_{2n-3,m-p})\}
$$
$$
\hat\tau^{n,m}_{2n+k}=-(e^{n-1,m}_{2n+k}+\beta_n e^{n-1,m}_{2n+k-1}+\alpha_n e^{n-2,m}_{2n+k-2})
+\sum^{m}_{p=0}b^{(n-1)}_{k+1,p}c_{2n-1,m-p} 
$$
$$
+b^{(n-1)}_{k,p}c_{2n,m-p}+\beta_n(b^{(n-1)}_{k,p}c_{2n-1,m-p}
+b^{(n-1)}_{k-1,p}c_{2n,m-p})
$$
$$
+\alpha_n(b^{(n-2)}_{k+1,p}c_{2n-3,m-p} +b^{(n-2)}_{k,p}c_{2n-2,m-p} +b^{(n-2)}_{k-1,p}c_{2n-1,m-p}+
b^{(n-2)}_{k-2,p}c_{2n,m-p})
$$
for $k=1,...,n$.
Next expand the remaining terms using the definitions in (\ref{checkAnm}), (\ref{checkBnm}) 
and the identities in (\ref{recidab}):
$$
\check A_{n,m}-\check B_{n,1}C_{2n,m-1}-...--\check B_{n,m}C_{2n,0}
=\{x\check\beta^{n,m}_0A_{n-1,0}+x^2\check\beta^{n,m}_0A_{n-2,0}\}
$$
$$-
\{x\check\beta^{n,m}_0B_{n-1,0}+x^2\check\beta^{n,m}_0B_{n-2,0}\}C_{2n,0}
+x\{\sum^{m-1}_{p=1}\check\beta^{n,m}_p(E_{n-1,p}
$$
$$
-\sum^p_{s=0}B_{n-1,s}(c_{2n-1,p-s}x^{2n-1}+c_{2n,p-s}x^{2n}))\}
+x^2\{\sum^{m-1}_{p=1}\check\alpha^{n,m}_p(E_{n-2,p}
$$
$$-\sum^p_{s=0}B_{n-2,s}(c_{2n-3,p-s}x^{2n-3}
+c_{2n-2,p-s}x^{2n-2}+c_{2n-1,p-s}x^{2n-1}+c_{2n,p-s}x^{2n}))\}
$$
$$
=
\{x\check\beta^{n,m}_0A_{n-1,0}+x^2\check\alpha^{n,m}_0A_{n-2,0}\}
-\{x\check\beta^{n,m}_0B_{n-1,0}+x^2\check\alpha^{n,m}_0B_{n-2,0}\}C_{2n,0}
$$
$$
-(\check \sigma^{(m)}_{2n-1}x^{2n-1}+\check \sigma^{(m)}_{2n}x^{2n})+
\check\tau^{n,m}_{2n+1}x^{2n+1}+...+\check\tau^{n,m}_{3n}x^{3n}
$$
where
$$
\check \sigma^{(m)}_{2n-1}=\sum^{m-1}_{p=1}\check\alpha^{n,m}_{p}\{-e^{n-2,p}_{2n-3}
+\sum^p_{s=0}b^{(n-2)}_{0,s}c_{2n-3,p-s}\}
$$
\begin{eqnarray*}
\check\sigma^{(m)}_{2n}&=&\sum^{m-1}_{p=1}\check\beta^{n,m}_{p}\{-e^{n-1,p}_{2n-1}
+\sum^p_{s=0}b^{(n-2)}_{0,s}c_{2n-1,p-s}\}\\
&&+
\check\alpha^{n,m}_{p}\{-e^{n-2,p}_{2n-2}
+\sum^p_{s=0}(b^{(n-2)}_{0,s}c_{2n-2,p-s}+b^{(n-2)}_{1,s}c_{2n-3,p-s})\}
\end{eqnarray*}
$$
\check\tau^{n,m}_{2n+k}=\sum^{m-1}_{p=1}\check\beta^{n,m}_{p}\{-e^{n-1,p}_{2n+k-1}
+\sum^p_{s=0}( b^{(n-1)}_{k,s}c_{2n-1,p-s}+b^{(n-1)}_{k-1,s}c_{2n,p-s})\}
$$
$$+
\check\alpha^{n,m}_{p}\{-e^{n-2,p}_{2n+k-2}
+\sum^p_{s=0}(b^{(n-2)}_{k+1,s}c_{2n-3,p-s}+b^{(n-2)}_{k,s}c_{2n-2,p-s}
$$
$$+b^{(n-2)}_{k-1,s}c_{2n-1,p-s}+b^{(n-2)}_{k-2,s}c_{2n,p-s})\}
$$
Recombining both parts, we get:
$$
E_{n,m}=\{x\check\beta^{n,m}_0A_{n-1,0}+x^2\check\alpha^{n,m}_0A_{n-2,0}\}
-\{x\check\beta^{n,m}_0B_{n-1,0}+x^2\check\alpha^{n,m}_0B_{n-2,0}\}C_{2n,0}
$$
$$
-(\sigma^{(m)}_{2n-1}x^{2n-1}+\sigma^{(m)}_{2n}x^{2n})+
(\hat\tau^{n,m}_{2n+1}+\check\tau^{n,m}_{2n+1})x^{2n+1}+...+
(\hat\tau^{n,m}_{3n}+\check\tau^{n,m}_{3n})x^{3n}
$$
where $
\sigma^{(m)}_{2n-1}=\hat \sigma^{(m)}_{2n-1}+\check\sigma^{(m)}_{2n-1}, 
\sigma^{(m)}_{2n}=\hat\sigma^{(m)}_{2n}+\check\sigma^{(m)}_{2n}$.
As in \ref{appA}, we then write:
$$
\sigma^{(m)}_{2n-1}x^{2n-1}+\sigma^{(m)}_{2n}x^{2n}=A^{(n)}[\sigma^{(m)}_{2n-1},\sigma^{(m)}_{2n}]
-B^{(n)}[\sigma^{(m)}_{2n-1},\sigma^{(m)}_{2n}]C_{2n,0}+...
$$
Thus, considering the forms of $A^{(n)}[\sigma^{(m)}_{2n-1},\sigma^{(m)}_{2n}]$, $B^{(n)}[\sigma^{(m)}_{2n-1},\sigma^{(m)}_{2n}]$ 
given in (\ref{Ann}), (\ref{Bnn}), it is clear that
if we define $\check\beta^{n,m}_0$, $\check\alpha^{n,m}_0$ as:
\begin{eqnarray}
\check\beta^{n,m}_0
&=&F^{(n)}_{2n-1}\sigma^{(m)}_{2n-1}+F^{(n)}_{2n}\sigma^{(m)}_{2n},\label{checkbetnm}\\
\check\alpha^{n,m}_0
&=&F^{(n-1)}_{2(n-1)}\sigma^{(m)}_{2n-1},\label{checkalphanm}
\end{eqnarray}
for $n\geq 1$,
then (\ref{Enm}) follows with
$$
e^{n,m}_{2n+k}=\hat\tau^{n,m}_{2n+k}+\check\tau^{n,m}_{2n+k}
-\tau^{(n)}_{2n+k}[\sigma^{(m)}_{2n-1},\sigma^{(m)}_{2n}],\ k=1,...,n,
$$
thus confirming that the equations in (\ref{v8}) are satisfied.\qedwhite

For the induction over $n$ the coefficients in $\{A_{n,p}\}$, $\{B_{n,p}\}$ for $p=0,...,m$
can be stored as matrices and generated simultaneously according to the proposition. It is also useful to store the parameters
$\{\check\alpha^{k,p}_0\}$, $\{\check\beta^{k,p}_0\}$ as matrices (for $k=1,...,n$, $p=1,...,m$) to assist with the 
induction over $m$. For example, evaluating the parameters according to (\ref{checkbetnm}), (\ref{checkalphanm}) in the case $n=1$:
\begin{eqnarray}
\check\beta^{1,m}_0&=&\frac{-1}{c_{1,0}}\{c_{2,m}+c_{1,0}c_{1,m}
+\check\beta^{1,m}_1c_{1,1}+...+\check\beta^{1,m}_{m-1}c_{1,m-1}\},\label{checkbetnm1}\\
\check\alpha^{1,m}_0&=&\hat c^{(m)}_1=b_{0,0}c_{1,m}+...+b_{0,m}c_{1,0},\label{checkalphanm1}
\end{eqnarray}
one can use the rules $\check\beta^{1,m}_p=\check\beta^{1,m-p}_0$ for $p<m$. 
These parameters can also be used to evaluate $A_{1,m}=a_{0,m}+a^{(1)}_{1,m}x$, $B_{1,m}=b_{0,m}+b^{(1)}_{1,m}x$,
where
\begin{eqnarray}
a^{(1)}_{1,m}&=&\check\alpha^{1,m}_0+\check\beta^{1,m}_1a_{0,1}+...+\check\beta^{1,m}_{m-1}a_{0,m-1}+\beta_1a_{0,m}\label{a11m}\\
b^{(1)}_{1,m}&=&\check\beta^{1,m}_0+\check\beta^{1,m}_1b_{0,1}+...+\check\beta^{1,m}_{m-1}b_{0,m-1}+\beta_1b_{0,m}\label{b11m}
\end{eqnarray}
Given the form of $\check\beta^{1,m}_0$, $\check\alpha^{1,m}_0$ in (\ref{checkbetnm1}), (\ref{checkalphanm1}),
it is clear that both $a^{(1)}_{1,m}$ and $b^{(1)}_{1,m}$ will be linear in $b_{0,1},...,b_{0,m}$.
It is also obvious from the formula for $\check\sigma^{(m)}_{2n}$ given above that $\check\beta^{n,m}_0$ will depend linearly on
$b_{0,1},...,b_{0,m}$ for $n=2$ and nonlinearly for $n>2$ (and likewise $\check\alpha^{n,m}_0$ for $n>1$).
Though this can certainly complicate the evaluation $b_{0,s}\to b^{(m)}_{0,s}$ for $s=1,...,p$ in (\ref{evalAB}), 
the example in the next section demonstrates that this process can be simplified in particular cases.
\section{Application}
\label{s3}
One area in which multivariate Pad\'e approximants naturally arise is as approximate solutions
for ordinary differential equations. 
Consider the (singular) {\em Riccati differential equation}:
\begin{equation}
xw'-\beta w+\beta w^2 +\alpha x=0,\ w(0)=0\ (w(1)=0)
\label{riccati}
\end{equation}
where $\alpha$ and $\beta$ are positive constants. As a singular first order differential equation,
according to a theorem of Malmquist (see Hille \cite{hille}), the {\em general solution} 
(for non-integer $\beta$) can be written as a bivariate series in $x$ and $x^\beta$, i.e.
\begin{equation}
w(x)=f(x,x^\beta)=\sum_{n,m=0}^\infty c_{n,m}x^nx^{m\beta};\ c_{0,0}=0
\label{riccbiv}
\end{equation}
where $c_{0,1}$ is an arbitrary parameter, specified by the condition $w(1)=0$.
However, there is a different approach to solving (\ref{riccati}) which suggests that 
this bivariate series is not the most natural form of the solution.
Using the substitution $u(x)=w'(x)/w(x)$ one derives
a second order linear equation for $u(x)$ which can be solved exactly in terms of Bessel functions, in 
turn giving the solution for $w(x)$ as a ratio of Bessel functions, i.e.
\begin{equation}
w(x)=\frac{z}{2\beta}\left\{
\frac{J_{\beta-1}(z)-C Y_{\beta - 1}(z)}{J_{\beta}(z)-C Y_{\beta}(z)}
\right\}
\label{riccati_soln}
\end{equation}
where
\[
z=2(\alpha\beta x)^{1/2}\ \ \hbox{and}\ \ C=
\frac{J_{\beta-1}(2\sqrt{\alpha\beta})}{Y_{\beta-1}(2\sqrt{\alpha\beta})}
\]
Expanding the Bessel functions as series reveals that the solution in (\ref{riccati_soln}) can be written as a 
rational function of the form:
\begin{equation}
w=\frac{\sum^\infty_{m=1} a_{m,0} x^m+x^\beta(\sum^\infty_{m=0} a_{m,1}x^m)}{
1+\sum^\infty_{m=1} b_{m,0} x^m+x^\beta(\sum^\infty_{m=0} b_{m,1}x^m)}
\label{riccrat}
\end{equation}
where
\begin{eqnarray*}
a_{m,0}&=&\frac{(-1)^m(\alpha\beta)^m}{(m-1)!(-\beta)_{m+1}},\ \ b_{m,0}=\frac{(-1)^m(\alpha\beta)^m}{m!(1-\beta)_m},
\\
a_{m,1}&=&\frac{\pi(-1)^{m}(1/\Gamma(\beta)+C\Gamma(1-\beta)\cos((\beta-1)\pi)/\pi)(\alpha\beta)^{m+\beta}}{
\beta C\Gamma(\beta)(\beta)_{m}m!}\\
b_{m,1}&=&\frac{\pi(-1)^{m}(1/\Gamma(\beta+1)+C\Gamma(-\beta)\cos(\beta\pi)/\pi)(\alpha\beta)^{m+\beta}}{
C\Gamma(\beta)(\beta+1)_{m}m!}
\end{eqnarray*}
and $(c)_n$ is the Pochhammer symbol.

The form of the solution in (\ref{riccrat}) suggests the use of multivariate Pad\'e approximants
as approximate solutions,
rather than truncating the bivariate series in (\ref{riccbiv}). In particular,  the 
{\em left}- 
and {\em right}-$(n,1)$ Pad\'e approximants
$$w^L_{n,1}(x)=f^L_{n,1}(x,x^\beta),\ \ w^R_{n,1}(x)=f^R_{n,1}(x,x^\beta),$$
both seem natural choices in
that they both share the same form, i.e. from (\ref{lpade}) and (\ref{rpade}):
$$\frac{A_{n,0}+A_{n,1}x^\beta}{B_{n,0}+B_{n,1}x^\beta},$$
and thus might be expected to converge to the exact solution as $n\to\infty$. 
In the following we first calculate each of the $w_{n,0}$, $w^L_{n,1}$, and $w^R_{n,1}$ using the (General)
algorithms from sections \ref{s1} and \ref{s2} that would be applicable to any bivariate series. We then consider 
how each of the methods can be refined in this particular case (Riccati)
to improve the speed and efficiency of the algorithm. 

The univariate part $w_{n,0}=A_{n,0}/B_{n,0}$ is the same for both $w^L_{n,1}$ and $w^R_{n,1}$, 
but there are several different ways of generating the polynomials $A_{n,0}$, $B_{n,0}$. {\em Algorithm 1 (General)}:
Calculate the parameters $\alpha_n$ from (\ref{alpn}), $\beta_n$ from  (\ref{betn}),
and errors $e^{(n)}_{2n+k}$ from (\ref{En0}); 
simplifying each using the recursive formulas in (\ref{jaca})-(\ref{jacb}). {\em Algorithm 2 (Riccati)}: Use the
explicit formulas for $\alpha_n$ and $\beta_n$ which are available in this case, i.e.
$$\alpha_n=\frac{-\alpha^2\beta^2}{(\beta-2n+1)(\beta-2n+2)^2(\beta-2n+3)},\ 
\beta_n=\frac{-2\alpha\beta}{(\beta-2n)(\beta-2n+2)}
$$
for $n \geq 2$, to simplify the recursion formulas. {\em Algorithm 3 (Riccati)}: The final method is to generate the coefficients
$b^{(n)}_{k,0}$ (and similarly $a^{(n)}_{k,0}$) directly from previous values, i.e.
\begin{equation}
b^{(n)}_{k,0}=P^{(0)}_{n,k} b^{(n)}_{k-1,0},\ 
a^{(n)}_{k,0}=Q^{(0)}_{n,k} a^{(n)}_{k-1,0}
\label{bnk0iter}
\end{equation}
where
$$
P^{(0)}_{n,k}=\frac{-\alpha\beta(2n-2k+1)(n-k+1)}{k(n-(k-1)/2)(\beta-(2n-k+1))(\beta-k)},
$$
$$
Q^{(0)}_{n,k}=\frac{-\alpha\beta(2n-2k+3)(n-k+1)}{(k-1)(n-(k-1)/2)
(\beta-(2n-k+2))(\beta-k)}
$$
These can be derived from the explicit formulas for $\beta_n$ $\alpha_n$ above 
and the recursive formulas in (\ref{jaca}), (\ref{jacb}). This is the most efficient approach since not only does
it generate the coefficients directly (without the induction), the formulas do not
require simplification.

There are several options for evaluating the {\em left}-Pad\'e approximant 
$w^L_{n,1}$, reflecting those for the univariate part above. 
{\em Algorithm 1 (General)}: 
This method starts by using Algoritm 1 above for the univariate part, $w_{n,0}$. 
The parameters $\check\alpha^{n,1}_0$, $\check\beta^{n,1}_0$ and errors $e^{n,1}_{2n+k}$
are then calculated using the formulas in section \ref{s2}, and from these the polynomials
$A_{n,1}$, $B_{n,1}$ are constructed from (\ref{Anm}), (\ref{Bnm}) for $m=1$.
{\em Algorithm 2 (Riccati)}:
This method starts by using Algorithm 2 above for the univariate part, $w_{n,0}$.
Then use the explicit formulas which are available in this case 
(and explicit formulas for $\beta_n$, $\alpha_n$ above), i.e. for $n\geq 2$:
$$
\check\beta^{n,1}_0=\frac{\beta(1-2\beta)_{2n-3}[2(2n-1)^2+\beta(-8(n-1)-3+2\beta)]\beta_nc_{0,1}}{
n(2n-1)!}
$$
$$
\check\alpha^{n,1}_0=\frac{8(n-2+\beta)(1-2\beta)_{2n-5}[8(n-1)^2+\beta(-8(n-1)-1+2\beta)]\alpha_nc_{0,1}}{(2n-1)!}
$$

For the {\em right}-Pad\'e approximant $w^R_{n,1}$ one calculates the parameters $\check\beta^{k,1}_0$, $\check\alpha^{k,1}_0$
using (\ref{checkbetnm1}) and (\ref{checkalphanm1}) respectively.
{\em Algorithm 1 (General)}: This method starts by using Algorithm 2 above for the univariate part, $w_{n,0}$. 
Since the coefficients $a^{(1)}_{k,1}$, $b^{(1)}_{k,1}$ in (\ref{a11m}), (\ref{b11m}) are linear in $b_{1,0},...,b_{k,0}$,
one can use the recursive formulas
for $a^{(n)}_{k,0}$, $b^{(n)}_{k,0}$ to derive similar formulas:
\begin{eqnarray}
a^{(n,1)}_{k,1}&=&a^{(n-1,1)}_{k,1}+\beta_na^{(n-1,1)}_{k-1,1}+\alpha_na^{(n-2,1)}_{k-2,1},\label{ank1induct}\\
b^{(n,1)}_{k,1}&=&b^{(n-1,1)}_{k,1}+\beta_nb^{(n-1,1)}_{k-1,1}+\alpha_nb^{(n-2,1)}_{k-2,1},\label{bnk1induct}
\end{eqnarray}
for $k<n$.
In fact, the formulas in (\ref{ank1induct}), (\ref{bnk1induct}) also work for $k\geq n$ using the obvious definitions:
$$
a^{(n,1)}_{n+s,1}=\{c_{n+s,1}+c_{n+s-1,1}b^{(n)}_{1,0}+...+c_{s,1}b^{(n)}_{n,0}
\}  +\check\beta^{n+s,1}_1a^{(n)}_{1,0}+...
+\check\beta^{n+s,1}_{n}a^{(n)}_{n,0}
$$
$$
b^{(n,1)}_{n+s,1}=\check\beta^{n+s,1}_0  +\check\beta^{n+s,1}_1b^{(n)}_{1,0}+...
+\check\beta^{n+s,1}_{n}b^{(n)}_{n,0}
$$
kor $s>0$.
{\em Algorithm 2 (Riccati)}: This method starts by using Algorithm 3 above for the univariate part, $w_{n,0}$.
Using the explicit formulas for $\beta_n$, $\alpha_n$ and the inductive approach in (\ref{bnk1induct}) 
one can derive iterative relations for $b^{(n,1)}_{k,1}$ (and $a^{(n,1)}_{k,1}$) which are similar to
those for $b^{(n)}_{k,0}$, $a^{(n)}_{k,0}$ in (\ref{bnk0iter}), i.e.
\begin{equation}
b^{(n,1)}_{k,1}=P^{(1)}_{n,k} b^{(n,1)}_{k-1,1},\ 
a^{(n,1)}_{k,1}=Q^{(1)}_{n,k} a^{(n,1)}_{k-1,1}
\label{bnk1iter}
\end{equation}
where
$$
P^{(1)}_{n,k}=\frac{-\alpha\beta (\beta-(n-(k-1)))(2\beta-(2n-2k+1))}{
k(\beta-(n-(k-1)/2))
(\beta-(2n-k+1))(k+\beta)}
$$
$$
Q^{(1)}_{n,k}=\frac{-\alpha\beta (\beta-(n-k+1))(2\beta-(2n-2k+3))}{k(\beta-(n-(k-2)/2))
(\beta-(2n-(k-1)))(k+\beta)}
$$
Thus, all of the coefficients in $B_{n,1}$ can be generated from $b_{0,1}$ without
using the induction. Note that the formulas in (\ref{bnk1iter}) also hold for $k\geq n$, allowing one to use
$b^{(n-1,1)}_{n,1}=P^{(1)}_{n-1,n}b^{(n-1,1)}_{n-1,1}$ to improve the induction in Algorithm 1.

To represent the error in the approximants
as a simple scalar value we first calculate the estimates of $c_{0,1}$ from the equations
$$A_{n,0}(1)+A^{L/R}_{n,1}(1)=0,$$
the solutions of which are denoted by $c^{L_n}_{0,1}$ and $c^{R_n}_{0,1}$, respectively.
Thus we represent the error as the absolute difference between these values and the exact value, $c^E_{0,1}$
(the $a_{0,1}$ coefficient from (\ref{riccrat})),
as given in Table~\ref{tab1}.
\begin{table}
\centering
\begin{tabular}{|c|c|c|}
\hline
$n$ & $|c^{L_n}_{0,1}-c^E_{0,1}|$ & $|c^{R_n}_{0,1}-c^E_{0,1}|$ \\
\hline
$1$ & $7.34\ (.045,.002)$ & $7.34\ (.002, .004)$ \\ 
$2$ & $1.68\ (.008,.002)$ & $1.53\ (.002,.012)$ \\
$3$ & $1.53\times 10^{-3}\ (.012,.005)$ & $1.83\times 10^{-2}\ (.8,.016)$ \\
$4$ & $2.45\times 10^{-3}\ (11,.008)$ & $3.61\times 10^{-5}\ (2.,.029)$ \\
$5$ & $5.51\times 10^{-4}\ (54.7,.008)$ & $9.69\times 10^{-9}\ (4.6,.073)$ \\ 
$6$ & $1.76\times 10^{-4}\ (?,.472)$ & $3.26\times 10^{-11}\ (11.87,.063)$ \\
$7$ & $7.00\times 10^{-5}\ (?,1.11)$ & $8.19\times 10^{-14}\ (39.2,.08)$\\
$8$ & $3.21\times 10^{-5}\ (?,1.7)$ & $2.10\times 10^{-16}\ (?,.11)$\\
$9$ & $1.63\times 10^{-5} (?,3.08)$ & $4.73\times 10^{-19}\ (?,.17)$\\
$10$ & $9.00\times 10^{-6}\ (?,4.88)$ & $9.37\times 10^{-22}\ (?,.21)$\\
\hline
\end{tabular} 
\caption{Errors in estimates of $c^E_{0,1}$.}
\label{tab1}
\end{table}
The time taken (up to a maximum allowable time of 1 minute) for each calculation 
using the general (Algorithm 1) and refined (Algorithm 2)
methods are given in brackets beside each error. From this it is obvious that the refined methods
are significantly faster than the general methods, and that the $w^R_{n,1}$ calculations
are faster than those for $w^L_{n,1}$.
The reasons for these differences are algebraic identities that were utilised for the refined methods
and the greater use of the Mathematica function
{\em Simplify} in the general methods to reduce the rational expressions.
It is also obvious from the errors themselves that $w^R_{n,1}(x)$ converges more quickly and uniformly
to the exact solution than $w^L_{n,1}(x)$.
This is explained by the fact (provable by induction) that $w^R_{n,m}=w^R_{n,1},\ m\geq 1$,
thus confining the errors in $w^R_{n,1}$ to higher orders.

\section{Conclusion}
As the Riccati example in the previous section illustrates, there are many different ways to
define and evaluate multivariate Pad\'e approximants and one must be prudent in considering the alternatives. 
For that example this meant comparing the {\em left-} and {\em right-}
Pad\'e approximants, $w^L_{n,1}$ and $w^R_{n,1}$; the latter proving to be a better approximation in terms of 
accuracy and speed of computation.
This type of analysis can easily be generalised to Riccati equations with polynomial
coefficients as in Fair, \cite{fair}, and further to Riccati
equations of second order (also known as the {\em Duffing Equation}), as in Fair, Luke, \cite{fairluke}.
These and other applications such as
the approximation of special functions (e.g. see Lutterodt, \cite{lutt} for
the Appell function)
should be solved in the complex domain (as in Hille, \cite{hille}) and allow for
{\em branch points} in solutions.
Similar methods based on the univariate Jacobi formula may also be applied to other definitions of 
multivariate Pad\'e approximants, such
as the diagonal approximants considered in Chisolm, \cite{chisolm} and Levin, \cite{levin}.
One could also compare the recursive method to others, 
including the QD-like algorithm in Cuyt, \cite{cuyt2}
and {\em Toeplitz} matrix methods like {\em Block-Levinson} by writing the equations 
in (\ref{v8}) as a linear system for the coefficients.
It would also be interesting to understand the recursive method in terms of the other 
approaches given in Cuyt, \cite{cuyt3}, particularly that of
{\em multivariate continued fractions} where the recursions seem to play a similar role 
to that in the univariate theory.
The bivariate Jacobi formulas in (\ref{Anm}), (\ref{Bnm}) can also be applied to
specific ordinary differential equations using the $\tau$-method, such as the 
Michaelis-Menten equation in Hegarty, \cite{hegarty1}, and this does have some 
advantages over the series-based approach here as one can use obvious simplification
patterns, similar to the univariate case in Hegarty, \cite{hegarty}.

\section*{Acknowledgments}

I wish to thank Dr Gerrard Liddell and Dr Stephen Taylor for their comments and suggestions.
This work originated during a Post-Doctoral Fellowship in the School of Pharmacy
at the University of Otago under the supervision of Professor Stephen Duffull.

\appendix
\section{Remainder Terms}\label{appA}
Let $\sigma_{2n-1}$, $\sigma_{2n}$ be arbitrary scalars and consider the canonical problem from section \ref{s2} 
of how to write the remainder terms in the form:
\begin{eqnarray}
\sigma_{2n-1}x^{2n-1}+\sigma_{2n}x^{2n}&=&A^{(n)}[\sigma_{2n-1},\sigma_{2n}]-B^{(n)}[\sigma_{2n-1},\sigma_{2n}]C_{2n,0}\label{n7b}\\
&&+
\tau^{(n)}_{2n+1}x^{2n+1}+...+\tau^{(n)}_{3n}x^{3n}
\nonumber
\end{eqnarray}
where $A^{(n)}[\sigma_{2n-1},\sigma_{2n}], B^{(n)}[\sigma_{2n-1},\sigma_{2n}]$ are polynomials of order $n$, 
both with zero constant components. These can be calculated recursively from the univariate Pad\'e polynomials, i.e.
\begin{eqnarray}
A^{(n)}[\sigma_{2n-1},\sigma_{2n}]
&=&xb^{{(n)}^*}_1 A_{n-1,0}+x^2b^{{(n-1)}^*}_1A_{n-2,0},\label{Ann}\\
B^{(n)}[\sigma_{2n-1},\sigma_{2n}]
&=&xb^{{(n)}^*}_1 B_{n-1,0}+x^2b^{{(n-1)}^*}_1B_{n-2,0},
\label{Bnn}
\end{eqnarray}
where the parameters $b^{{(n)}^*}_1$, $b^{{(n-1)}^*}_1$ and the errors $\tau^{(n)}_{2n+k}$ can be written in component form:
\begin{eqnarray*}
b^{{(n)}^*}_1&=&b^{{(n)}^*}_1[\sigma_{2n-1},\sigma_{2n}]=F^{(n)}_{2n-1}\sigma_{2n-1}+F^{(n)}_{2n}\sigma_{2n}\\
b^{{(n-1)}^*}_1&=&b^{{(n-1)}^*}_1[0,\sigma_{2n-1}]=F^{(n-1)}_{2(n-1)}\sigma_{2n-1}\\
\tau^{(n)}_{2n+k}&=&\tau^{(n)}_{2n+k}[\sigma_{2n-1},\sigma_{2n}]=\tau^{(n)}_{2n+k,2n-1}\sigma_{2n-1}+
\tau^{(n)}_{2n+k,2n}\sigma_{2n},
\end{eqnarray*}
for $\ k=1,...,n$, and where the compoments $F^{(n)}_{2n-1}$, $F^{(n)}_{2n}$, $\tau^{(n)}_{2n+k,2n-1}$,
$\tau^{(n)}_{2n+k,2n}$ are independent of $\sigma_{2n-1}$, $\sigma_{2n}$.
This starts with the case $n=1$, where
$$
F^{(0)}_0=1,\ F^{(1)}_1=0,\ F^{(1)}_2=-1/c_{1,0},\ \tau^{(1)}_{3,1}=0,\ \tau^{(1)}_{3,2}=-c_{2,0}/c_{1,0}
$$
Substituting the formulas from (\ref{Ann}), (\ref{Bnn}) back into (\ref{n7b}) and expanding orders $x^{2n-1}$, 
$x^{2n}$ gives the components:
$$
F^{(n)}_{2n-1}=\frac{-\tau^{(n-1)}_{2(n-1)+1,2(n-1)}}
{\{c_{2n-1,0}-e^{(n-1)}_{2(n-1)+1}\}},\ 
F^{(n)}_{2n}=-1/
\{c_{2n-1,0}-e^{(n-1)}_{2(n-1)+1}\}
$$
Similarly, expanding orders $x^{2n+k}$ gives error terms 
$\tau^{(n)}_{2n+k}=\tau^{(n)}_{2n+k}[\sigma_{2n-1},\sigma_{2n}]$
for $k=1,...,n$ in component form, i.e.
$$
\tau^{(n)}_{2n+k}=\{F^{(n)}_{2n-1}[b^{(n-1)}_{k-1,0}c_{2n,0}+b^{(n-1)}_{k,0}c_{2n-1,0}-e^{(n-1)}_{2n+k-1}]
+F^{(n-1)}_{2(n-1)}[b^{(n-2)}_{k-2,0}c_{2n,0}
$$
$$+b^{(n-2)}_{k-1,0}c_{2n-1,0}
+b^{(n-2)}_{k,0}c_{2n-2,0}+b^{(n-2)}_{k+1,0}c_{2n-3,0}-e^{(n-2)}_{2n+k-1}]\}\sigma_{2n-1}
$$
$$
+\{F^{(n)}_{2n}[b^{(n-1)}_{k-1,0}c_{2n,0}+b^{(n-1)}_{k,0}c_{2n-1,0}-e^{(n-1)}_{2n+k-1}]\}\sigma_{2n}
$$
$$
=\tau^{(n)}_{2n+k,2n-1}\sigma_{2n-1}+
\tau^{(n)}_{2n+k,2n}\sigma_{2n},
$$

\end{document}